\documentclass{article}

\usepackage{amsfonts}
\usepackage{amssymb}
\usepackage{amsthm}
\usepackage{amsmath}
\usepackage{newlfont}
\usepackage{graphicx}
\usepackage[all,arc]{xy}
\usepackage{epsfig}

\usepackage{multicol}
\usepackage[T1]{fontenc}

\newtheorem {proposition}{Proposition}[section]
\newtheorem {theorem}{Theorem}[section]
\newtheorem {lemma}{Lemma}[section]

\author{{\DJ{}or\dj{}e Barali\' {c}}\\ {\small Mathematical Institute SASA}\\[-2mm] {\small Belgrade, Serbia} \and Branko Grbi\'{c}\\ {\small Mathematical Grammar School}\\[-2mm] {\small Belgrade, Serbia}
\\[-2mm]\and \DJ{}or\dj{}e \v{Z}ikeli\'{c}\\ {\small Mathematical Grammar School}\\[-2mm] {\small Belgrade, Serbia}}

\title{\textit{Cinderella}, Quadrilaterals and Conics}
\date{}
\begin{document}
\maketitle

\begin{abstract} We study quadrilaterals inscribed and
circumscribed about conics. Our research is guided by experiments
in software \textit{Cinderella}. We extend the known results in
projective geometry of conics and show how modern mathematical
software brings new ideas in pure and applied mathematics.
Poncelet theorem for quadrilaterals is proved by elementary means
together with Poncelet's grid property.

\end{abstract}

\renewcommand{\thefootnote}{}
\footnotetext{This research was supported by the Grant 174020 of
the Ministry for Education and Science of the Republic of Serbia
and Project Math Alive of the Center for Promotion of Science,
Serbia and Mathematical Institute SASA.}

\section{Introduction}

This paper is one in the serial of our forthcoming papers in
geometry of curves, combinatorics and dynamic systems. The use of
software \textit{Cinderella} is common for all of them and our aim
is to show that the good software is more than a box with nice
examples and calculations. The smart use could lead us not only to
discovering new results, but it gives the complete and correct
proofs! In this sense \textit{Cinderella} could go beyond the
limits of geometry of conics and mechanical experiments, even to
the curves of higher degree and abstract combinatorics, geometry
and topology.

The positive experience with \textit{Cinderella} in the paper
\textit{Illumination of Pascal's Hexagrammum and Octagrammum
Mysticum} by Barali\'{c} and Spasojevi\'{c}, \cite{BarIgi}
encouraged us to continue the research. The problems we study are
strongly influenced by very inspirative paper \textit{Curves in
Cages: an Algebro-geometric Zoo} of Gabriel Katz printed in
American Mathematical Monthly, \cite{Katz}. Many important
questions in dynamical systems and combinatorics have their
equivalents in the terms of algebraic curves. Richard Schwartz and
Serge Tabachnikov in \cite{Taba} asked for the proof of Theorem
4.c. They found the theorem studying the pentagram maps,
introduced in \cite{Schva}. This is still open hypothesis and
could be reformulated in the question about curves.

We have not found the proof for Schwartz and Tabachnikov Theorem
4.c but during recent work we discovered new interesting facts
about quadrilaterals inscribed and quadrilaterals circumscribed
about conic. Theorems about quadrilaterals and conics are usually
known like degenerate cases of Pascal and Brianchon theorems. In
\cite{BarIgi} Barali\'{c} and Spasojevi\'{c} proved some new
results about two quadrilaterals inscribed in a conic. However, in
this paper we study more complicated structures involving both
tangents at the vertices and the side lines of quadrilateral. We
start from the degenerate form of Pascal and Brianchon theorems
for the quadrilateral and then we discover new interesting points,
conics and loci.

The objects are studied by elementary means. Some of the results
are in particular the corollary of Great Poncelet Theorem for the
case when $n$-gon is quadrilateral. Here we give the short proof
for this case. Some special facts about this special case are
explained as well.

Finally, we compare two theorems - Mystic Octagon theorem for the
case of two quadrilaterals and Poncelet Theorem for the
quadrilaterals. Both of them have in common that they state that
certain $8$ points coming from two quadrilaterals inscribed in a
conic lie on the same conic. While the first one is pure
algebro-geometric fact, the latter involves much deeper structure
of the space and can not be seen naturally as the special case of
the first. Thus, we could not find 'Theorem of all theorems for
conics in projective geometry' and elementary surprises in
projective geometry like those in \cite{Taba} could come as the
special case of different general statements.

\section{From Pascal to Brocard Theorem}

In this section we show how Pascal theorem for hexagon (1639)
inscribed in a conic degenerates to Brocard theorem for the
quadrilateral inscribed in a circle. All results here are well
known and are part of the standard olympiad problem solving
curriculum, but our aim is to illustrate the power of degeneracy
tool and prepare the background for the next sections.

\begin{figure}[h!h!]
\centerline{\epsfig{figure=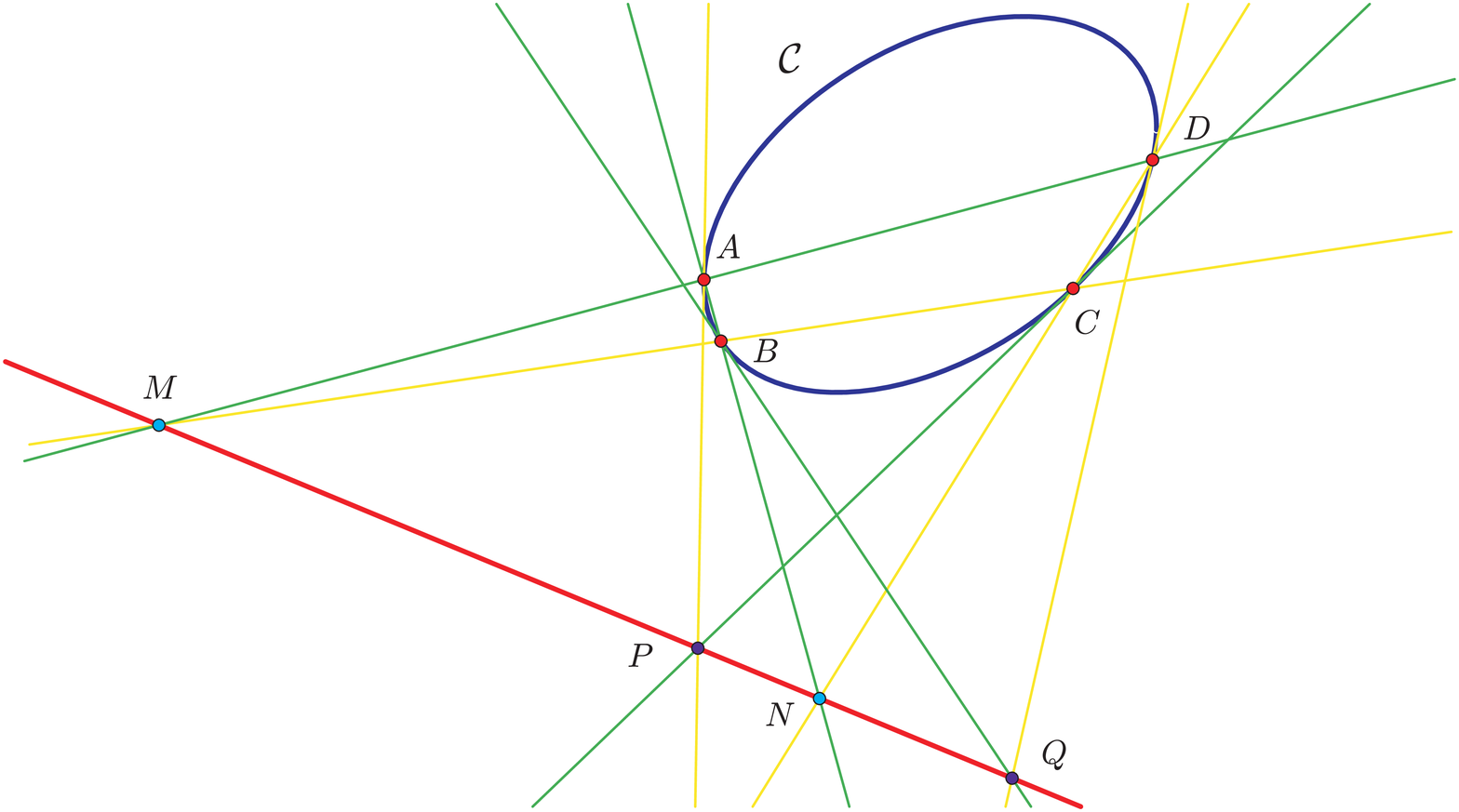,width=0.6\textwidth}}
\caption{Lemma \ref{slikacetvorougao}} \label{cet1}
\end{figure}

\begin{lemma}\label{slikacetvorougao} Let $A B C D$ be a quadrilateral inscribed in a conic $\mathcal{C}$ and let $M$ be the intersection point of the lines $A D$ and $B C$, $N$
be the intersection point of the lines $A B$ and $C D$, $P$ be the
intersection point of the tangents to $\mathcal{C}$ at $A$ and
$C$, and $Q$ be the intersection point of the tangents to
$\mathcal{C}$ at $B$ and $D$. Then, the points $M$, $N$, $P$ and
$Q$ are collinear (see Figure \ref{cet1}).
\end{lemma}

\noindent {\bf Proof:} Apply Pascal theorem to degenerate hexagon
$A A B C C D$ and we get the points $M$, $N$ and $P$ are
collinear. Apply Pascal theorem to degenerate hexagon $A B B C D$
and we get the points $M$, $N$ and $Q$ are collinear. \hfill
$\square$

\medskip

\begin{figure}[h!h!]\vspace*{-0.5 cm}
\centerline{\epsfig{figure=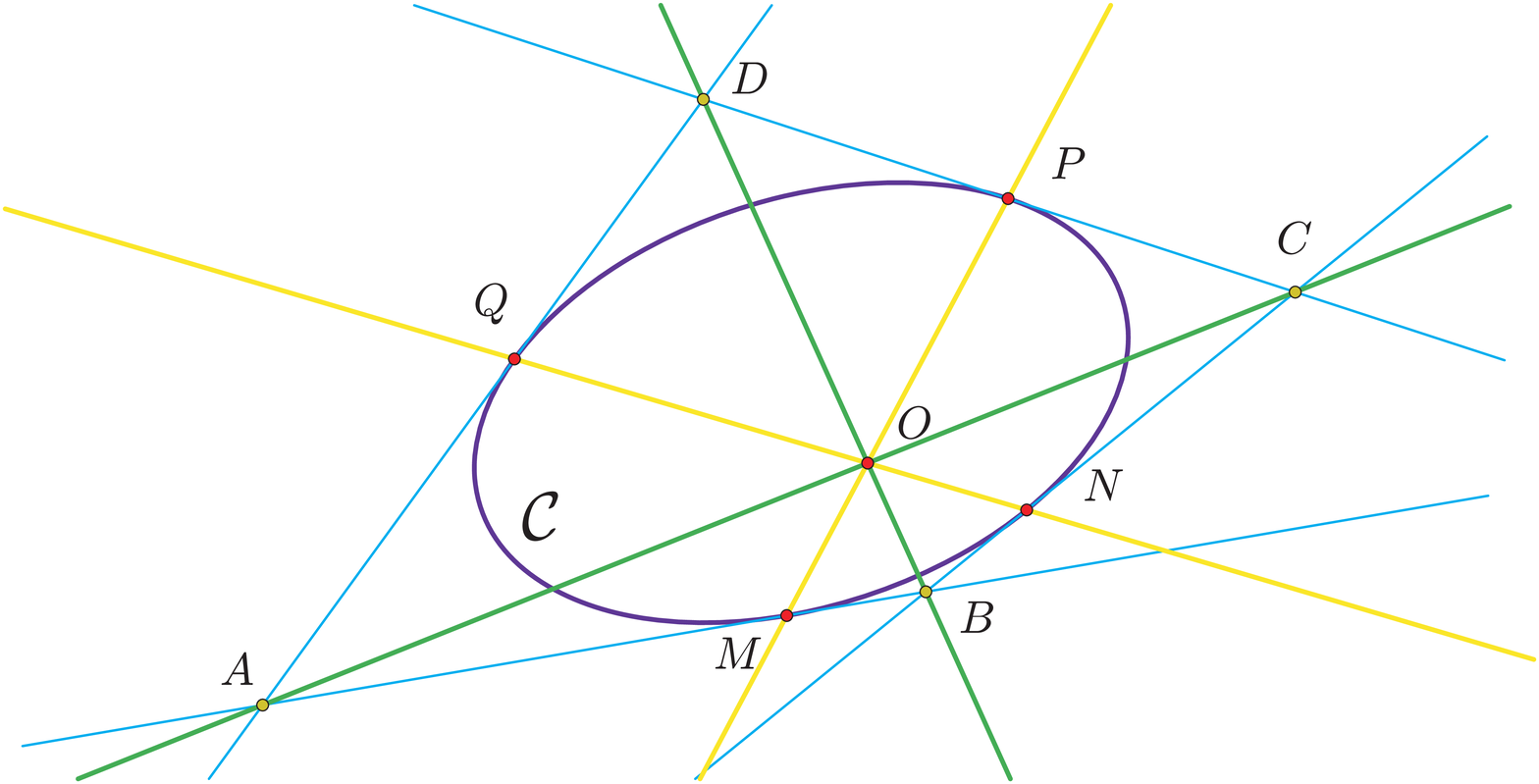,width=0.5\textwidth}}
\caption{Lemma \ref{cetb}} \label{cetbslika}
\end{figure}

Dual statement to Lemma \ref{slikacetvorougao} is the following:

\begin{lemma}\label{cetb} Let conic $\mathcal{C}$ touch the sides $A B$, $B C$, $C D$ and $D A$ of a quadrilateral $A B C D$ in the points $M$, $N$, $P$ and $Q$, respectively.
Then the lines $A C$, $B D$, $M P$ and $N Q$ pass through the same
point $O$ (see Figure \ref{cetbslika}).
\end{lemma}

The Lemmas \ref{slikacetvorougao} and \ref{cetb} will be used to
prove the other interesting relations among the lines and points
that naturally occur in a quadrilateral inscribed in conics
configurations. Many points are going to be introduced so we are
going to organize labels of the points.

Let $A_1 A_2 A_3 A_4$ be a quadrilateral inscribed in a conic
$\mathcal{C}$ and let $M_1$ be the intersection point of the lines
$A_1 A_2$ and $A_3 A_4$, $M_2$ of $A_2 A_3$ and $A_4 A_1$ and
$M_3$ of $A_3 A_1$ and $A_2 A_4$. Let $N_3$ be the intersection
point of the tangent lines to the conic at $A_1$ and $A_3$, $P_3$
of the tangents at $A_2$ and $A_4$, $N_2$ of the tangents at $A_1$
and $A_4$, $P_2$ of the tangents at $A_2$ and $A_3$, $N_1$ of the
tangents at $A_1$ and $A_2$ and $P_1$ of the tangents at $A_3$ and
$A_4$. Let $U_1$ and $U_2$ be the points where tangents from $M_1$
touch $\mathcal{C}$, and analogously $V_1$, $V_2$ and $W_1$, $W_2$
for the points $M_2$ and $M_3$ respectively.

Lemma \ref{slikacetvorougao} states that the points $M_1$, $M_2$,
$N_3$, $P_3$ are collinear, and also the points $M_2$, $M_3$,
$N_1$, $P_1$ and $M_3$, $M_1$, $N_2$ and $P_2$. Denote this three
lines by $m_3$, $m_1$ and $m_2$, respectively. We are going to
prove that $U_1$ and $U_2$ lie on the line $m_3$, $V_1$ and $V_2$
on the line $m_2$ and $W_1$ and $W_2$ on $m_1$ - so that $m_1$,
$m_2$ and $m_3$ are the polar lines of the points $M_1$, $M_2$ and
$M_3$ with respect to $\mathcal{C}$.

\begin{lemma}\label{satang} The points $U_1$, $M_2$, $U_2$ and
$M_3$ are collinear.
\end{lemma}

\begin{figure}[h!h!]
\centerline{\epsfig{figure=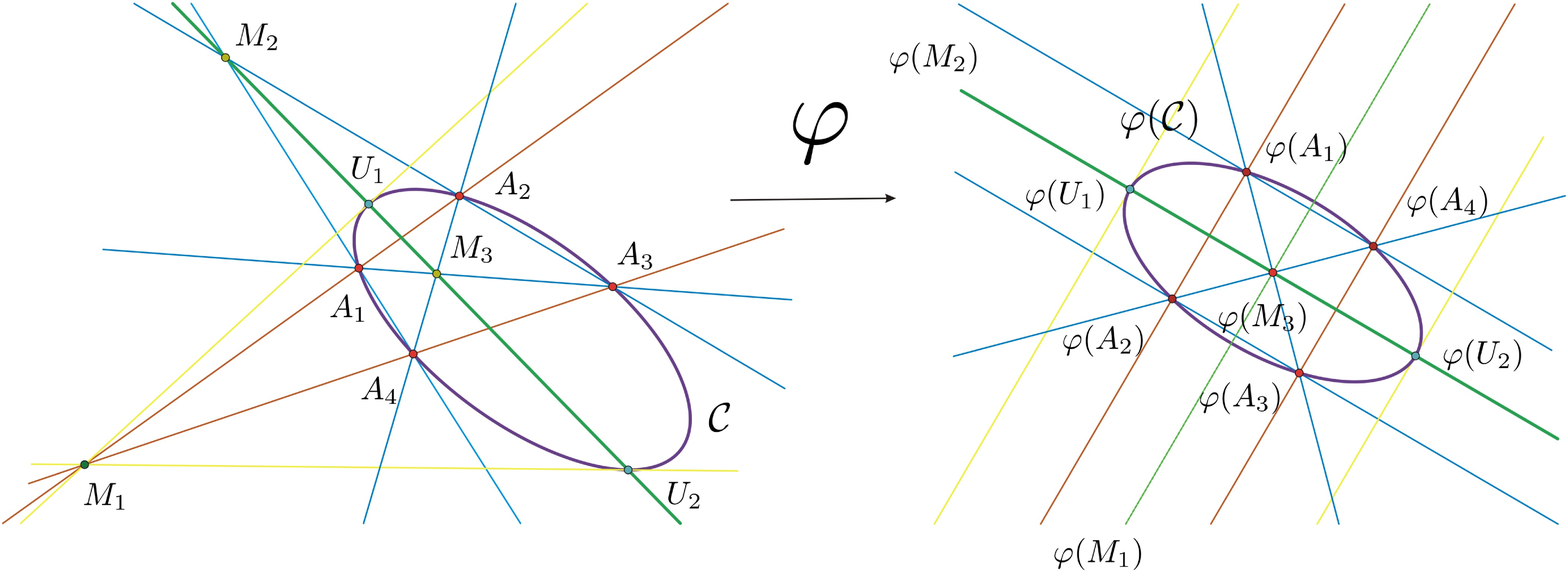,width=\textwidth}}
\caption{Lemma \ref{satang}} \label{zatan}
\end{figure}

\noindent {\bf Proof:} There is a projective transformation
$\varphi$ that maps the points $A_1$, $A_2$, $A_3$ and $A_4$ onto
the vertices of a square. Thus $\varphi(M_3)$ is the center of a
square with vertices $\varphi (A_1)$, $ \varphi (A_2)$, $ \varphi
(A_3)$ and $ \varphi (A_4)$. The points $\varphi(M_1)$ and
$\varphi (M_2)$ are at infinity. There is a unique way to inscribe
the square into the conic, and the lines $\varphi(A_1)
\varphi(A_2)$ and $\varphi(A_1) \varphi(A_4)$ are parallel to the
axes of the conic $\varphi (\mathcal{C})$. The points $\varphi
(U_1)$ and $\varphi (U_2)$ must be mapped onto the axis parallel
to the line $\varphi(A_1) \varphi(A_4)$.

Now the points $\varphi (U_1)$, $\varphi (U_2)$, $\varphi (M_2)$
and $\varphi (M_3)$ lie on the axis of conic $\varphi
(\mathcal{C})$. Consequently, the points $U_1$, $M_2$, $U_2$ and
$M_3$ then lie at the same line. \hfill $\square$

\begin{figure}[h!h!]
\centerline{\epsfig{figure=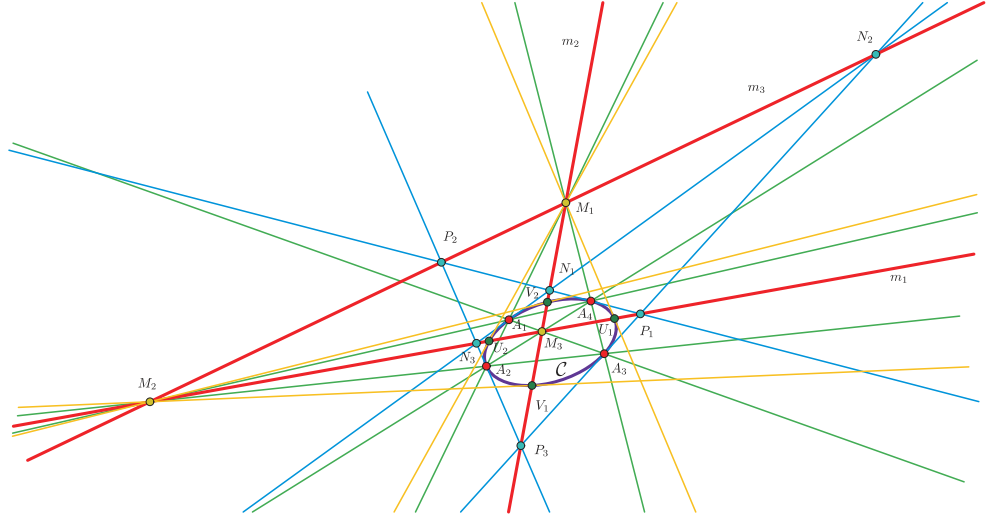,width=\textwidth}}
\caption{Quadrilateral inscribed in a conic} \label{triprave}
\end{figure}

Lemma \ref{satang} clearly implies the analogous statement for the
lines $m_2$ and $m_3$. This is the classical theorem of the
projective geometry and a very useful tool. Some other facts about
a quadrilateral inscribed in conics are going to be proved in the
next sections.

\begin{figure}[h!h!]
\centerline{\epsfig{figure=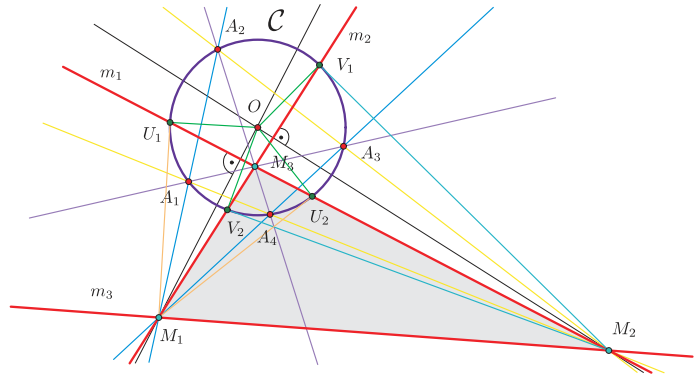,width=0.7\textwidth}}
\caption{Brocard theorem} \label{tripravek}
\end{figure}

In the end, we treat one very special case - when the conic
$\mathcal{C}$ is a circle. Projective geometry gives us the plenty
of techniques. For example, in the proof of Lemma \ref{satang} we
used the projective transformation. We have already described
degeneracy tool when we take some limit cases of polygons
inscribed (or circumscribed) in a conic. It is good to keep in
mind that conic could degenerate itself for example to the two
lines. This is a way to get interesting configurations of points
and lines.

The configuration \ref{triprave} in the case of a circle has nice
a property which is known as the Brocard theorem. Let $O$ be the
center of a circle $\mathcal{C}$. Then the quadrilateral $M_1 U_1
O U_2$ is deltoid and we get $M_1 O\perp m_1$. Similarly, $M_2
O\perp m_2$. Thus:

\begin{theorem}[Brocard theorem] Let $O$ be the center of circumscribed circle of a
cyclic quadrilateral $A_1 A_2 A_3 A_4$. Then $O$ is the
orthocenter of triangle $\triangle M_1 M_2 M_3$.
\end{theorem}

\section{More lines and pencils of lines}

We continue in the same manner. The lines and the pencils of lines
we study came from various degenerations of the vertices of
hexagon inscribed in a conic. Let us remind that configuration
associated with 60 Pascal lines has been described in \cite{Ladd},
\cite{Vero} and \cite{BarIgi}. All results from this section could
be obtained as the certain degenerate case. But we are going to
treat them by elementary means.

Let $T_1$ be the point of intersection of the line $A_3 A_4$ and
tangent at $A_1$ to $\mathcal{C}$, $T_2$ of $A_4 A_1$ and tangent
at $A_2$, $T_3$ of $A_1 A_2$ and tangent at $A_3$ and $T_4$ of
$A_2 A_3$ and tangent at $A_4$. Let $X_1$ be the point of
intersection of the line $A_2 A_3$ and tangent at $A_1$ to
$\mathcal{C}$, $X_2$ of $A_3 A_4$ and tangent at $A_2$, $X_3$ of
$A_4 A_1$ and tangent at $A_3$ and $X_4$ of $A_1 A_2$ and tangent
at $A_4$. Let $Y_1$ be the point of intersection of the line $A_2
A_3$ and tangent at $A_1$, $Y_2$ of $A_1 A_4$ and tangent at
$A_2$, $Y_3$ of $A_1 A_4$ and tangent at $A_3$ and $Y_4$ of $A_2
A_3$ and tangent at $A_4$.

\begin{proposition}\label{pro1} The following $16$ triples of points are
collinear: $(M_1, Y_1, Y_2)$, $(M_1, Y_3, Y_4)$, $(M_1, X_3,
T_4)$, $(M_1, X_1, T_2)$, $(M_2, Y_1, Y_4)$, $(M_2, Y_2, Y_3)$,
$(M_2, X_4, T_1)$, \\ $(M_2, X_2, T_3)$, $(M_3, T_1, T_3)$, $(M_3,
X_2, X_4)$, $(M_3, X_1, X_3)$, $(M_3, T_2, T_4)$, $(X_2, Y_3,
T_4)$, \\ $(X_1, Y_2, T_3)$, $(X_3, Y_4, T_1)$, $(X_4, Y_1, T_2)$.
\end{proposition}

\noindent {\bf Proof:} The collinearity of the points $M_1$, $X_3$
and $T_4$ follows from the Pascal theorem for degenerate hexagon
$A_1 A_4 A_4 A_3 A_3 A_2$, the collinearity of the points $M_1$,
$Y_3$ and $Y_4$ from degenerate hexagon $A_1 A_3 A_3 A_4 A_4 A_2$
and the collinearity of the points $X_2$, $Y_3$ and $T_4$ from
degenerate hexagon $A_2 A_3 A_3 A_4 A_4 A_2$. The proof for the
rest is analogous. \hfill $\square$

\begin{figure}[h!h!]
\centerline{\epsfig{figure=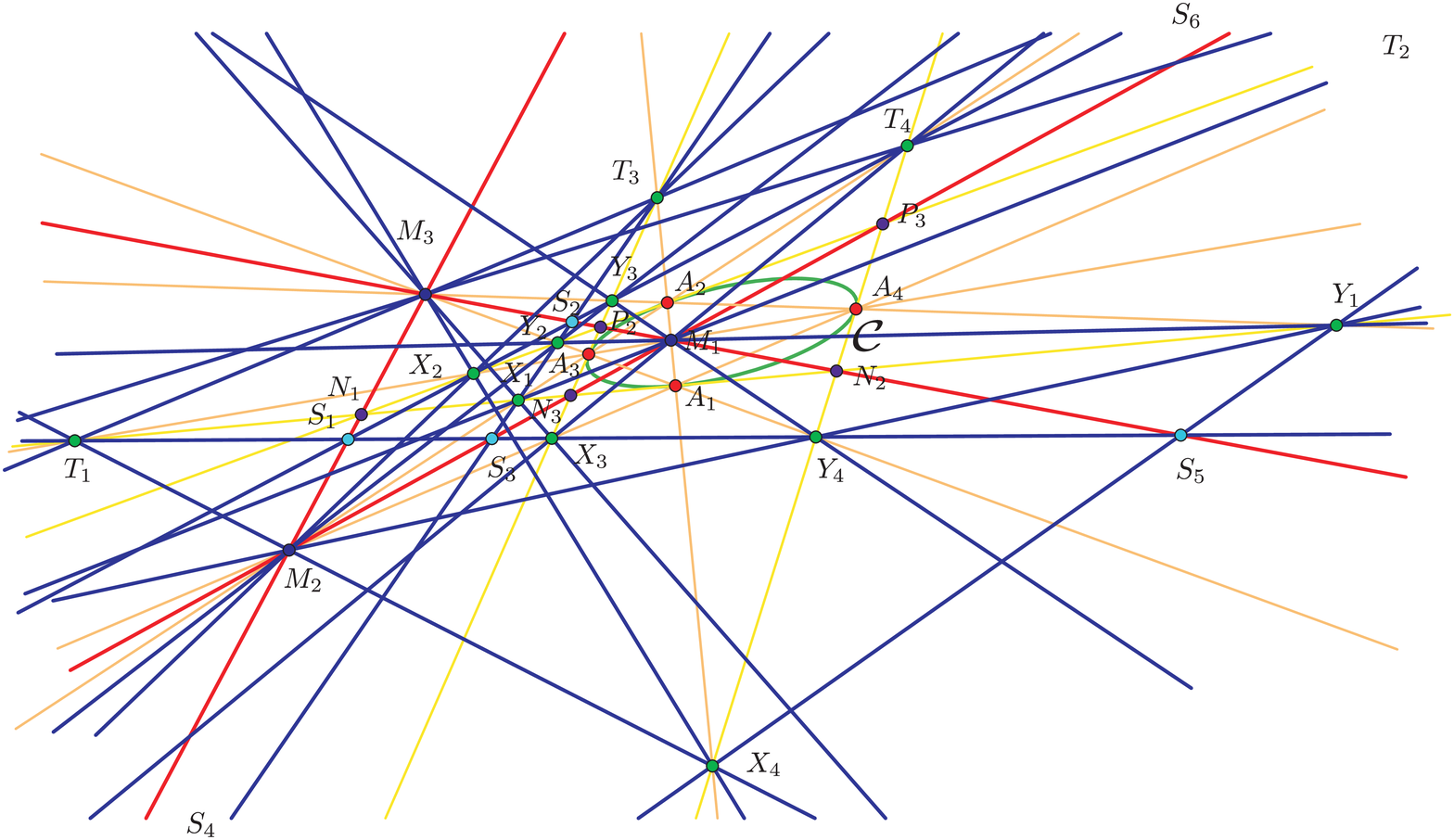,width=\textwidth}}
\caption{Propositions \ref{pro1} and \ref{pro2}} \label{slika1}
\end{figure}

\begin{proposition}\label{pro2} The following $6$ triples of lines are concurrent:\\ $(M_2 M_3, X_2 Y_3, X_3 Y_4 )$,
$(M_1 M_3, X_1 Y_2, X_2 Y_3)$, $(M_1 M_2, X_1 Y_2, X_3 Y_4)$,\\
$(M_2 M_3, X_1 Y_2, X_4 Y_1)$, $(M_1 M_3, X_4 Y_1, X_3 Y_4)$,
$(M_1 M_2, X_4 Y_1, X_2 Y_3)$.
\end{proposition}

\noindent {\bf Proof:} By Lemma \ref{oktt} (to be proved in the
next section) the points $X_1$, $X_2$, $X_3$, $X_4$, $T_1$, $T_2$,
$T_3$ and $T_4$ lie on the same conic. From Pascal theorem for the
hexagon $T_1 X_3 X_1 T_2 X_4 X_2$ we get that lines $M_1 M_3$,
$X_4 Y_1$ and $X_3 Y_4$ are concurrent. Analogously for other
triples. \hfill $\square$

\medskip

Define the points as the intersections of the lines: $B_1=l (A_2
V_1)\cap l (A_1 V_2)$, $C_1=l (A_1 V_1)\cap l (A_2 V_2)$, $D_1=l
(A_3 V_1)\cap l (A_4 V_2)$, $E_1=l (A_4 V_1)\cap l (A_3 V_2)$,
$B_3=l (A_4 V_1)\cap l (A_2 V_2)$, $C_3=l (A_4 V_2)\cap l (A_2
V_1)$, $D_3=l (A_1 V_1)\cap l (A_3 V_2)$, $E_3=l (A_1 V_2)\cap l
(A_3 V_1)$, $D_2=l (A_4 U_1)\cap l (A_1 U_2)$, $E_2=l (A_1
U_1)\cap l (A_4 U_2)$, $B_2=l (A_3 U_1)\cap l (A_2 U_2)$, $C_2=l
(A_2 U_1)\cap l (A_3 U_2)$, $F_3=l (A_4 U_1)\cap l (A_2 U_2)$,
$H_3=l (A_4 U_2)\cap l (A_2 U_1)$, $G_3=l (A_1 U_1)\cap l (A_3
U_2)$, $I_3=l (A_1 U_2)\cap l (A_3 U_1)$, $E_1=l (A_2 W_1)\cap l
(A_1 W_2)$, $F_1=l (A_1 W_1)\cap l (A_2 W_2)$, $G_1=l (A_3
W_1)\cap l (A_4 W_2)$, $H_1=l (A_4 W_1)\cap l (A_3 W_2)$, $H_2=l
(A_4 W_1)\cap l (A_1 W_2)$, $I_2=l (A_4 W_2)\cap l (A_1 W_1)$,
$F_2=l (A_2 W_1)\cap l (A_3 W_2)$ and $G_2=l (A_2 W_2)\cap l (A_3
W_1)$.

\begin{proposition}\label{superbijanson} The points $B_1$, $C_1$, $D_1$, $E_1$, $F_1$, $G_1$, $H_1$, $I_1$ lie on the line $M_2
M_3$. Similarly, the points $B_2$, $C_2$, $D_2$, $E_2$, $F_2$,
$G_2$, $H_2$, $I_2$ lie on the line $M_3 M_1$ and the points
$B_3$, $C_3$, $D_3$, $E_3$, $F_3$, $G_3$, $H_3$, $I_3$ lie on the
line $M_1 M_2$.
\end{proposition}

\begin{figure}[h!h!]
\centerline{\epsfig{figure=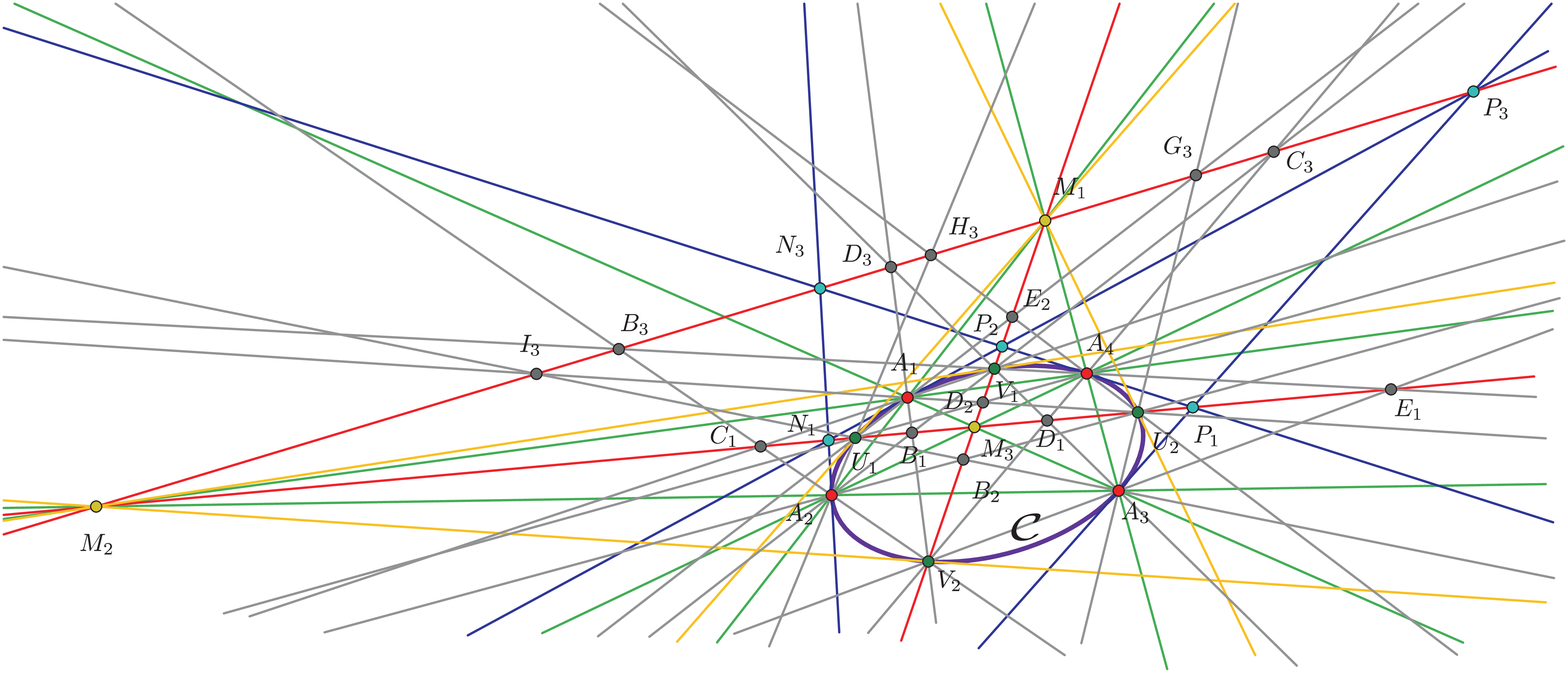,width=\textwidth}}
\caption{Propositions \ref{superbijanson}} \label{slika2}
\end{figure}

\noindent {\bf Proof:} Consider the quadrilateral formed by the
tangent lines to the conic $\mathcal{C}$ at the points $A_4$,
$A_2$, $V_1$ and $V_2$. Applying Lemma \ref{cetb}, we get that the
point $B_3$ lies on the line $M_1 M_2$. Analogously for other
points. \hfill $\square$

\section{Surprising conics}

In the upper sections many points were introduced. We have showed
some of them are collinear while some are the intersections of
certain lines. But some of them lie on the curves of degree two!

\begin{figure}[h!h!]
\centerline{\epsfig{figure=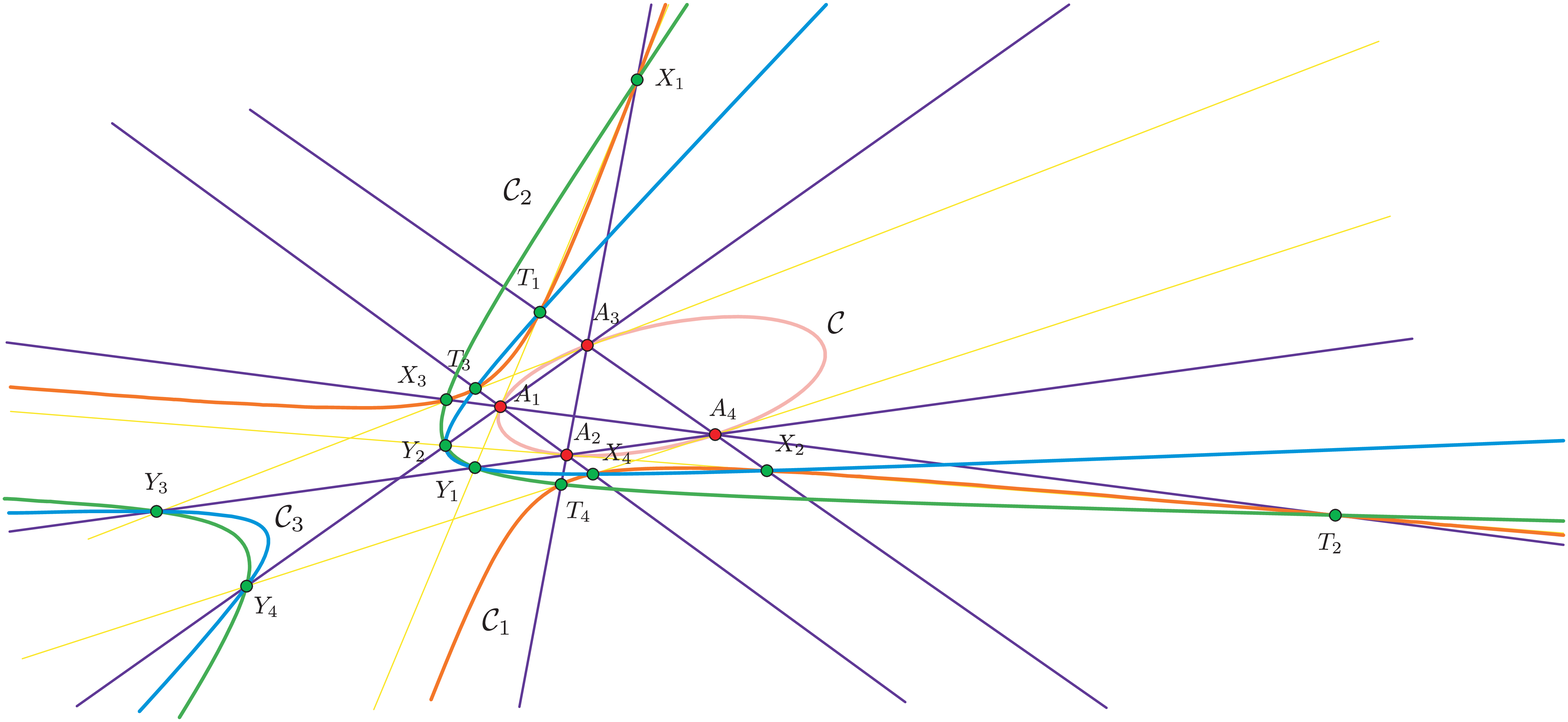,width=\textwidth}}
\caption{Propositions \ref{oktt}} \label{slika3}
\end{figure}

\begin{lemma}\label{oktt} The points $X_1$, $X_2$, $X_3$, $X_4$,
$T_1$, $T_2$, $T_3$ and $T_4$ lie on the same conic
$\mathcal{C}_1$; $Y_1$, $Y_2$, $Y_3$, $Y_4$, $X_1$, $X_3$, $T_2$,
and $T_4$ lie on the same conic $\mathcal{C}_2$; $T_1$, $T_3$,
$X_2$, $X_4$, $Y_1$, $Y_2$, $Y_3$ and $Y_4$ lie on the same conic
$\mathcal{C}_3$ (see Figure \ref{slika3}).
\end{lemma}

\noindent {\bf Proof:} This statement is the special case of the
Mystic Octagon theorem, the first time formulated by Wilkinson in
\cite{Wilk}. The first conic appears when we consider degenerate
octagon $A_1 A_2 A_2 A_3 A_3 A_4 A_4 A_1$, the second for $A_1 A_3
A_3 A_2 A_2 A_4 A_4 A_1$, and the third for $A_1 A_3 A_3 A_4 A_4
A_2 A_2 A_1$. \hfill $\square$

Let $J_{2i-1}$ be the intersection points of the tangents at
$X_{i-2}$ and $T_i$ on the conic $\mathcal{C}_1$, and $J_{2i}$ the
intersection points of the tangents at $X_{i-1}$ and $T_i$ (modulo
4), for $i=1$, $2$, $3$, $4$.  Then the following claim is true:

\begin{theorem}\label{prva} \begin{itemize}
    \item The lines $J_{i} J_{i+4}$, for
$i=1$, $2$, $3$, $4$ intersect at the point $M_3$.
    \item The lines $J_1 J_7$, $J_2 J_6$ and $J_3 J_5$ intersect at $M_1$
    and the lines $J_1 J_3$, $J_4 J_8$  and $J_5 J_7$ intersect at $M_2$.
    \item The lines $J_1 J_4$ and $J_2 J_5$ intersect at $A_1$,
    the lines $J_4 J_7$ and $J_3 J_6$ at $A_2$, the lines $J_6
    J_1$ and $J_5 J_8$ at $A_3$ and the lines $J_3 J_8$ and $J_2
    J_7$ at $A_4$.
    \item The intersection points $l (J_2 J_4)\cap l (J_6 J_8)$, $l (J_2 J_8)\cap l (J_4
    J_6)$, $l (J_3 J_6)\cap l (J_2 J_7)$, $l (J_5 J_8)\cap l (J_1
    J_4)$, $l (J_3 J_8)\cap l (J_4 J_7)$, $l (J_2 J_5)\cap l (J_1
    J_6)$ and $l (J_i J_{i+1})\cap l (J_{i+4} J_{i+5})$ for $i=1$, $2$, $3$,
    $4$ lie on the same line $M_1 M_2$.
    \item  The intersection points $l (J_4 J_5)\cap l (J_7 J_8)$ and $l (J_3 J_4)\cap l
    (J_1 J_8)$ lie on the same line $M_1 M_3$, the intersection points $l (J_2 J_3)\cap l (J_5 J_6)$ and $l (J_1 J_2)\cap l
    (J_6 J_7)$ lie on the same line $M_2 M_3$.
    \item The point $P_3$ lies on the line $J_3 J_7$ and the point
    $N_3$ on the line $J_1 J_5$.
    \item Three lines $J_{2i} J_{2i+4}$, $J_{2i+1} J_{2i-2}$ and $J_{2i-1}
    J_{2i+2}$ (modulo 8) are concurrent for $i=1$, $2$, $3$, $4$.
\end{itemize}
\end{theorem}

\begin{figure}[h!h!]
\centerline{\epsfig{figure=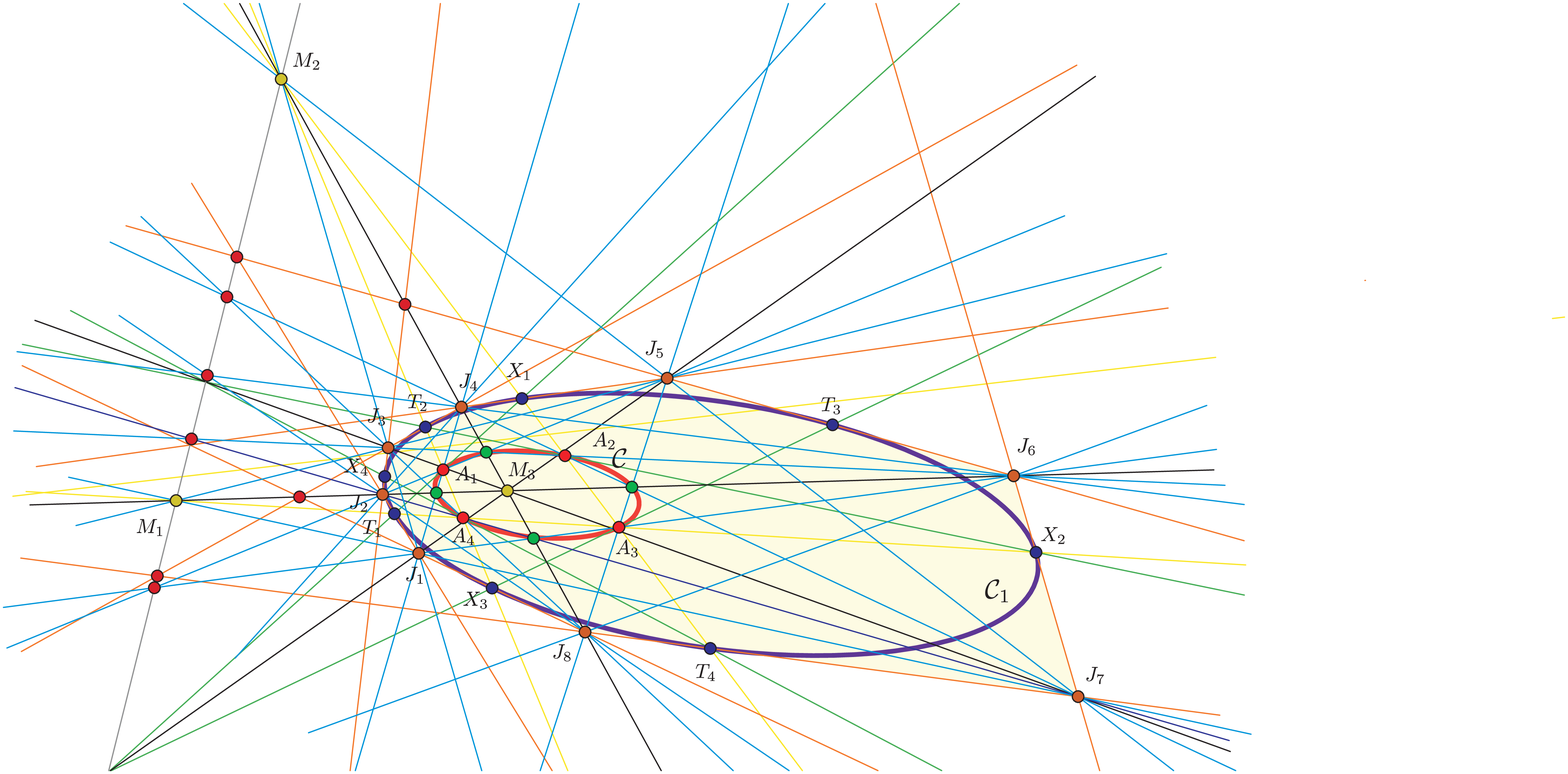,width=\textwidth}}
\caption{Theorem \ref{prva}} \label{slikat41}
\end{figure}

\noindent {\bf Proof:} Consider the quadrilateral formed by
tangents to $\mathcal{C}_1$ at $J_2$ and $J_6$. By Lemma
\ref{cetb} and Proposition \ref{pro1} the points $M_3$ and $M_2$
lie on the line $J_2 J_6$ (we could take the order of points
differently). Analogously, the lines $J_1 J_5$, $J_3 J_7$ and $J_4
J_8$ pass through the point $M_3$. In similar fashion we prove
other statements for the points $M_1$ and $M_2$, as well as the
points $N_3$ and $P_3$.

Lemma \ref{cetb} applied to the quadrilateral formed by tangents
to $\mathcal{C}_1$ at $J_2$ and $J_5$ proves that line $J_2 J_5$
pass through $A_1$. Similarly, $A_1$ belongs to the line $J_1
J_4$. Analogously, we prove the corresponding statements for the
points $A_2$, $A_3$ and $A_4$.

From Lemma \ref{slikacetvorougao} applied on the quadrilateral
$T_2 X_1 T_4 X_3$ and Proposition \ref{pro1}, it follows that the
intersection point of the lines $J_3 J_4$ and $J_7 J_8$ and the
intersection point of the lines $J_4 J_5$ and $J_8 J_1$ lie on the
line $M_1 M_2$. Then by Brianchon Theorem for the hexagon formed
by the tangents to $\mathcal{C}_1$ at $T_2$, $X_1$, $T_3$, $T_1$,
$X_3$ and $T_4$ the intersection point of the line $J_1 J_4$ and
$J_5 J_8$ lie on the line $M_1 M_2$. Analogously for the others.

Brianchon Theorem for the hexagon formed by the tangents to
$\mathcal{C}_1$ at $T_2$, $X_1$, $X_4$, $T_1$, $X_3$ and $T_4$
applies the concurrency of the lines $J_2 J_6$, $J_1 J_4$ and $J_5
J_8$. We use the similar argument for the rest of the proof.
\hfill $\square$

\medskip

\begin{figure}[h!h!]
\centerline{\epsfig{figure=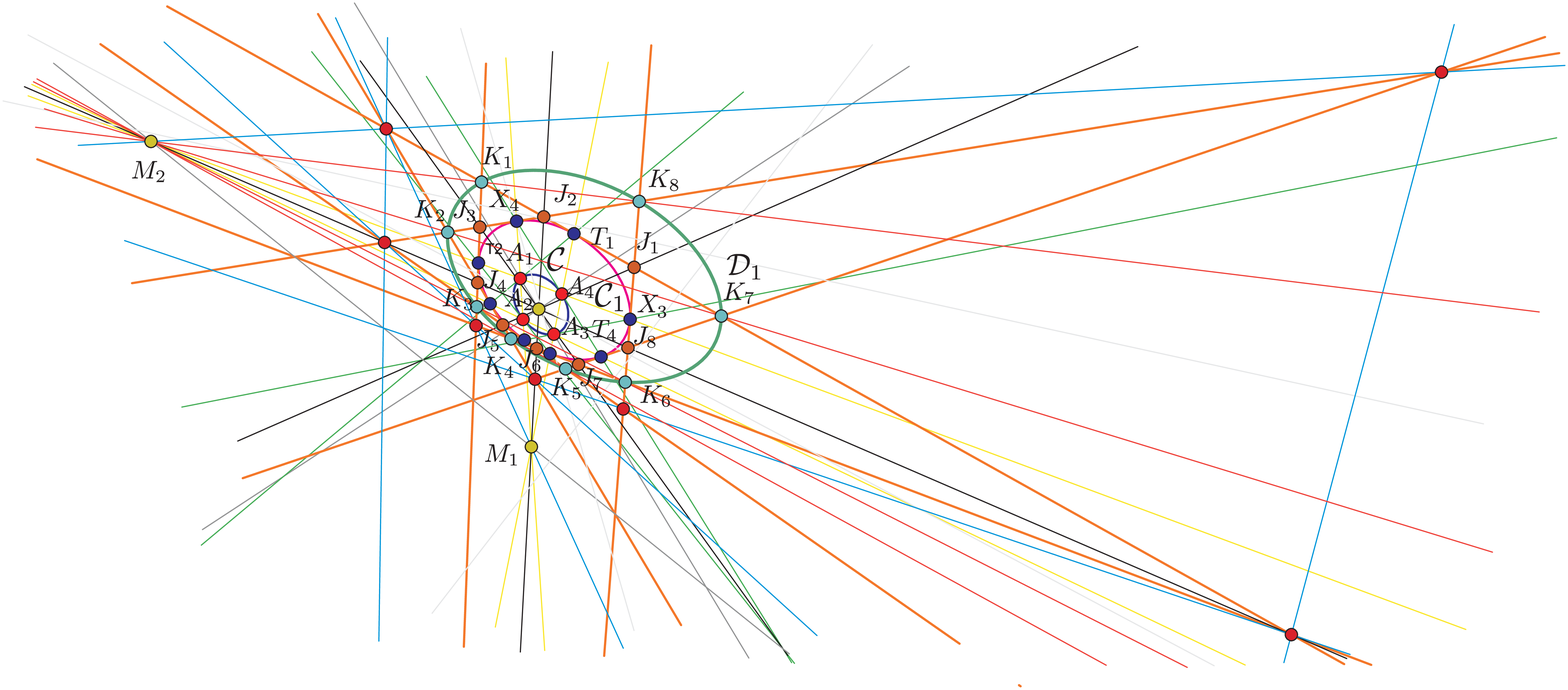,width=\textwidth}}
\caption{Theorem \ref{druga}} \label{slikat42}
\end{figure}

Let $K_i$ be the intersection points of lines $J_{i} J_{i+1}$ and
$J_{i+2} J_{i+3}$ (modulo 8) for $i=1$, $\dots$, $8$.

\begin{theorem}\label{druga} The points $K_i$ lie on the same conic
$\mathcal{D}_1$.
\end{theorem}

\noindent {\bf Proof:} It is not hard to prove that the lines $K_1
K_5$, $K_2 K_6$, $K_3 K_7$ and $K_4 K_8$ pass through the point
$M_3$, the lines $K_2 K_3$, $K_1 K_4$, $K_5 K_8$ and $K_6 K_7$
pass through the point $M_1$ and the lines $K_2 K_7$, $K_1 K_8$,
$K_3 K_6$ and $K_4 K_5$ pass through the point $M_2$. From the
collinearity of the points $M_1$, $J_2$ and $l (J_4 J_5)\cap l
(J_7 J_8)$ the points $K_1$, $K_2$, $K_4$, $K_5$, $K_7$ and $K_8$
lie on the same conic. Using the similar argument we show that
$K_2$, $K_4$, $K_5$, $K_6$, $K_7$ and $K_8$ lie on the same conic.
Because there is a unique conic determined by its 5 points then
all the points $K_1$, $K_2$, $K_4$, $K_5$, $K_6$, $K_7$ and $K_8$
are on the same conic. Then it is easy to prove that $K_3$ also
lies on the conic. \hfill $\square$

\medskip
Let $Z_1=l (M_1 U_1)\cap l (M_2 V_1)$, $Z_2=l (M_1 U_1)\cap l (M_2
V_2)$, $Z_3=l (M_1 U_2)\cap l (M_2 V_2)$ and $Z_1=l (M_1 U_2)\cap
l (M_2 V_1)$.

\begin{theorem}\label{tankon} The points $N_1$, $N_2$, $P_1$, $P_2$, $Z_1$,
$Z_2$, $Z_3$ and $Z_4$ lie on the same conic.
\end{theorem}

\begin{figure}[h!h!]
\centerline{\epsfig{figure=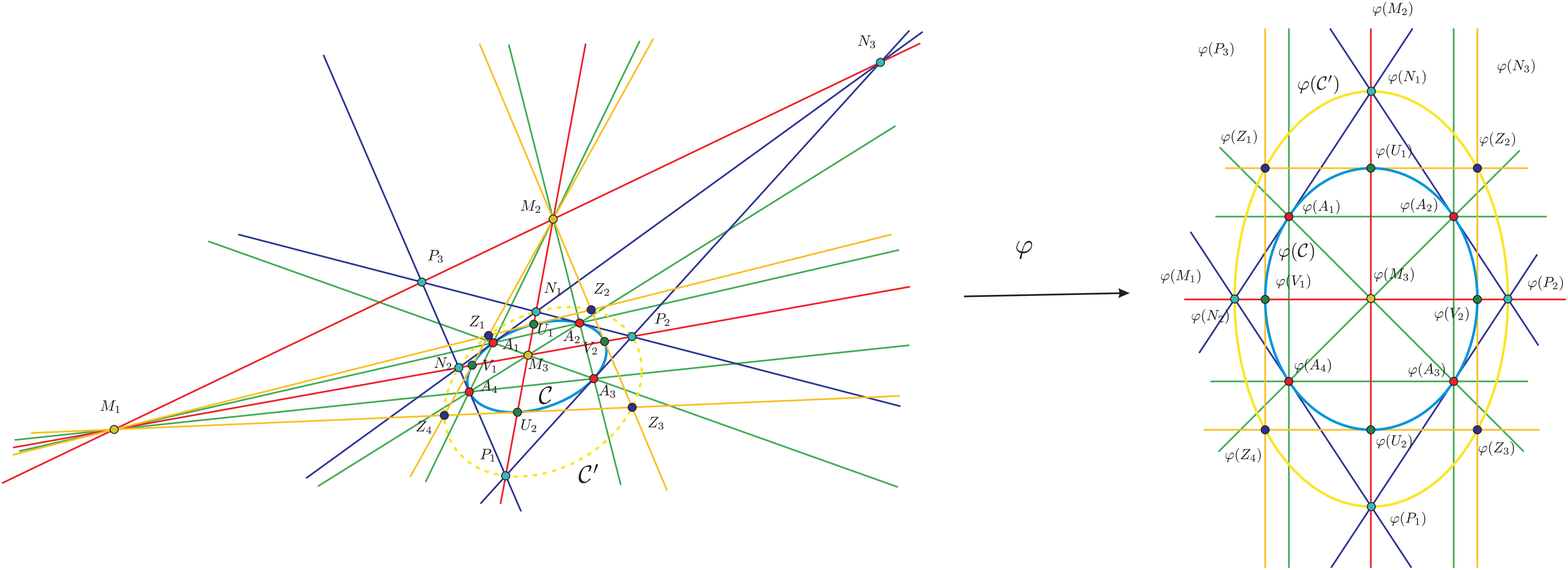,width=\textwidth}}
\caption{Theorem \ref{tankon}} \label{slikat43}
\end{figure}

\noindent {\bf Proof:} There exists projective transformation
$\varphi$ that maps vertices $A_1$, $A_2$, $A_3$ and $A_4$ onto
the vertices of a square. Then point $\varphi (M_3)$ is mapped
onto the center of conic $\varphi (\mathcal{C})$ and the lines
$\varphi (N_1) \varphi(P_1)$ and $\varphi (N_2) \varphi(P_2)$ are
the axes. The points $\varphi (U_1)$, $\varphi (U_2)$, $\varphi
(V_1)$ and $\varphi (V_2)$ also lie on the axes. As we could see
in Figure \ref{slikat43}, everything is symmetric and it is easy
to conclude that there is a conic through $\varphi(Z_1)$,
$\varphi(Z_2)$, $\varphi(Z_3)$, $\varphi(Z_4)$, $\varphi(N_1)$,
$\varphi(N_2)$, $\varphi(P_1)$ and $\varphi(P_2)$. \hfill
$\square$

Theorems \ref{prva}, \ref{druga} and \ref{tankon} associate new
conics to the quadrilateral inscribed in a conic. They have
interesting properties which will be explained in the next
section.

\section{Poncelet's quadrilateral porism}

Jean-Victor Poncelet's famous \textit{Closure theorem} states that
if there exists one $n$-gon inscribed in conic $\mathcal{C}$ and
circumscribed about conic $\mathcal{D}$ then any point on
$\mathcal{C}$ is the vertex of some $n$-gon inscribed in conic
$\mathcal{C}$ and circumscribed about conic $\mathcal{D}$.
Poncelet published his theorem in \cite{Ponce}. However, this
result influenced mathematics until nowadays. In recent book
\cite{DraRad} by Dragovic and Radnovic there are several proofs of
Closure theorem, it's generalizations as well as it's relations
with elliptic functions theory. The proof is not elementary for
general $n$, although in the case $n=3$ elegant proof could be
found in almost every monograph in projective geometry, see
\cite{Pra}.

Theorems \ref{druga} and \ref{tankon} are the special cases of
Poncelet theorem for $n=4$. Actually, quadrilaterals and conics in
them have poristic property. We kept the spirit of elementarity
through our paper and our agenda was: Firstlu, we experiment in
Cinderella, after that the proof is recovered by elementary tools
(again directly guided by Cinderella's tools). In the same style
we continue and offer direct analytic proof of Poncelet theorem
for quadrilaterals without using differentials and elliptic
functions.

\begin{figure}[h!h!]
\centerline{\epsfig{figure=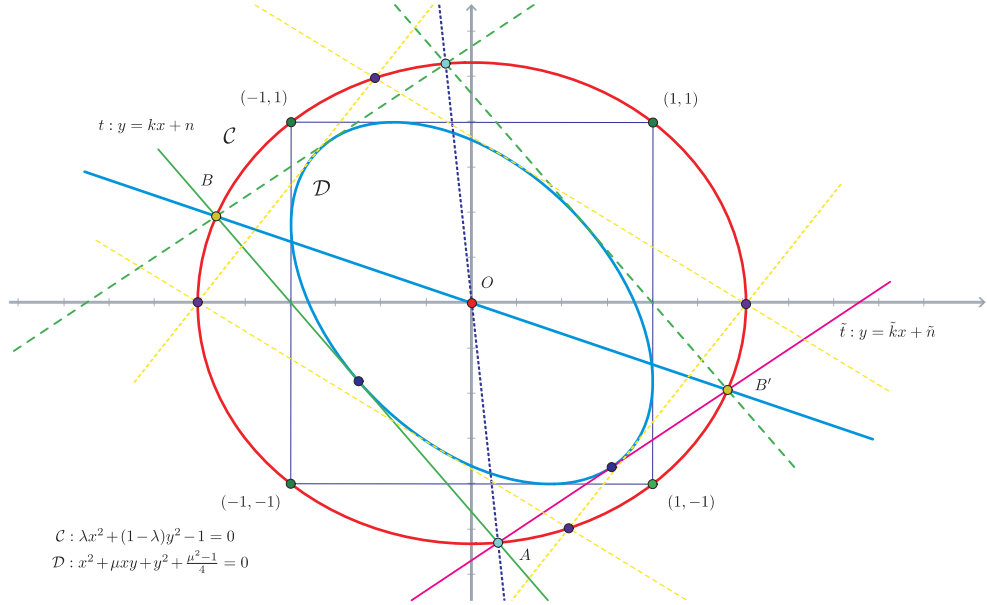,width=\textwidth}}
\caption{Lemma \ref{ponce} and Theorem \ref{Ponce4}}
\label{slikat44}
\end{figure}

\begin{lemma} \label{ponce} Let $\lambda$, $\mu$ be such that
conics $\mathcal{C} : \lambda x^2+(1-\lambda) y^2-1=0$ and
$\mathcal{D}: x^2+\mu x y+y^2+\frac{\mu^2-1}{4}=0$ are
non-degenerate. Let $A$ be a point on $\mathcal{C}$ and $B$ and
$B'$ be the intersections of the tangent lines from $A$ to
$\mathcal{D}$ with conic $\mathcal{C}$. Then the points $B$ and
$B'$ are symmetric with respect to the origin.
\end{lemma}

\noindent {\bf Proof:} Let a line $t: y=k x+n$ be a tangent line
to conic $\mathcal{D}$. The condition of tangency between $t$ and
$\mathcal{D}$ is \begin{equation}\label{c111} n^2 =k^2+m k+1.
\end{equation} The coordinates of the intersection points of $t$
and $\mathcal{C}$ are $$ (x_1,
y_1)=\left(\frac{-2(1-\lambda)kn-\sqrt{D}}{2(\lambda+(1-\lambda)k^2)},
k \cdot
\left(\frac{-2(1-\lambda)kn-\sqrt{D}}{2(\lambda+(1-\lambda)k^2)}\right)+n\right)$$
and $$ (x_2,
y_2)=\left(\frac{-2(1-\lambda)kn+\sqrt{D}}{2(\lambda+(1-\lambda)k^2)},
k \cdot
\left(\frac{-2(1-\lambda)kn+\sqrt{D}}{2(\lambda+(1-\lambda)k^2)}\right)+n\right),
$$ where $D=4 (\lambda-\lambda (1-\lambda) n^2+(1-\lambda) k^2)$.
It is necessary and enough to prove that the line through the
points $(-x_1, -y_1)$ and $(x_2, y_2)$ is tangent to
$\mathcal{D}$. This line has the equation $y=\tilde{k}
x+\tilde{n}$ where $\tilde{k}$ and $\tilde{n}$ could be calculated
as \begin{equation}\label{c113}
\tilde{k}=\frac{-\lambda}{(1-\lambda)k} \mbox{\, and \,}
 \tilde{n}=\frac{\sqrt{D}}{2 k (1-\lambda)}. \end{equation} We need
to check if $$\tilde{n}^2 =\tilde{k}^2+m \tilde{k}+1.$$ It is
directly verified that condition (\ref{c111}) multiplied by
$\frac{\lambda (1-\lambda)}{k^2 (1-\lambda)^2}$ finishes our
proof. \hfill $\square$

\begin{theorem}\label{Ponce4} Let $\mathcal{C}$ and $\mathcal{D}$
be conics such that there exists one quadrilateral inscribed in a
conic $\mathcal{C}$ and circumscribed about a conic $\mathcal{D}$.
Then any point on $\mathcal{C}$ is the vertex of some
quadrilateral inscribed in conic $\mathcal{C}$ and circumscribed
about conic $\mathcal{D}$.
\end{theorem}

\noindent {\bf Proof:} There exists a projective transformation
that maps the vertices of the quadrilateral inscribed in conic
$\mathcal{C}$ and circumscribed about conic $\mathcal{D}$ onto the
points $(1, 1)$, $(1, -1)$, $(-1, -1)$ and $(-1, 1)$ (in the
standard chart). Thus, conics $\mathcal{C}$ and $\mathcal{D}$ are
transformed in those with the equations as in Lemma \ref{ponce}.
Now the claim follows. \hfill $\square$

In fact, we proved more. All quadrilaterals with poristic property
with respect to $\mathcal{C}$ and $\mathcal{D}$ have the common
point of the intersection of diagonals (lines joining opposite
vertices) and the common line passing through the intersections of
opposite side lines. Our work in previous section, now could be
reviewed in the new light.

\begin{figure}[h!h!]
\centerline{\epsfig{figure=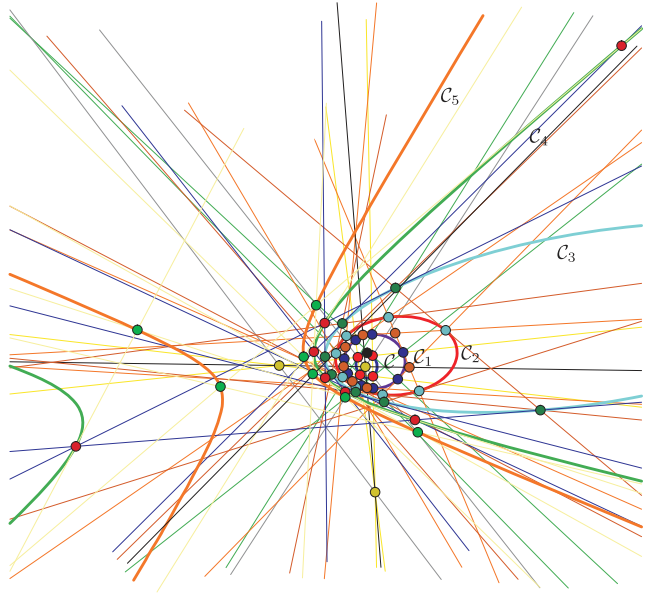, width=0.65\textwidth}}
\caption{The first five conics in the sequence} \label{slikat45}
\end{figure}

Theorems \ref{prva}, \ref{druga} i \ref{tankon} are obtained after
we defined certain points. If we apply the same procedure for
defining new points on the points and conics in Theorems, again we
come to similar conclusions. Thus, by repeating this procedure, we
obtain the infinite sequence of conics, see Figure \ref{slikat45}.
Every two consecutive conics in this sequence are Poncelet
4-connected.

\begin{figure}[h!h!]
\centerline{\epsfig{figure=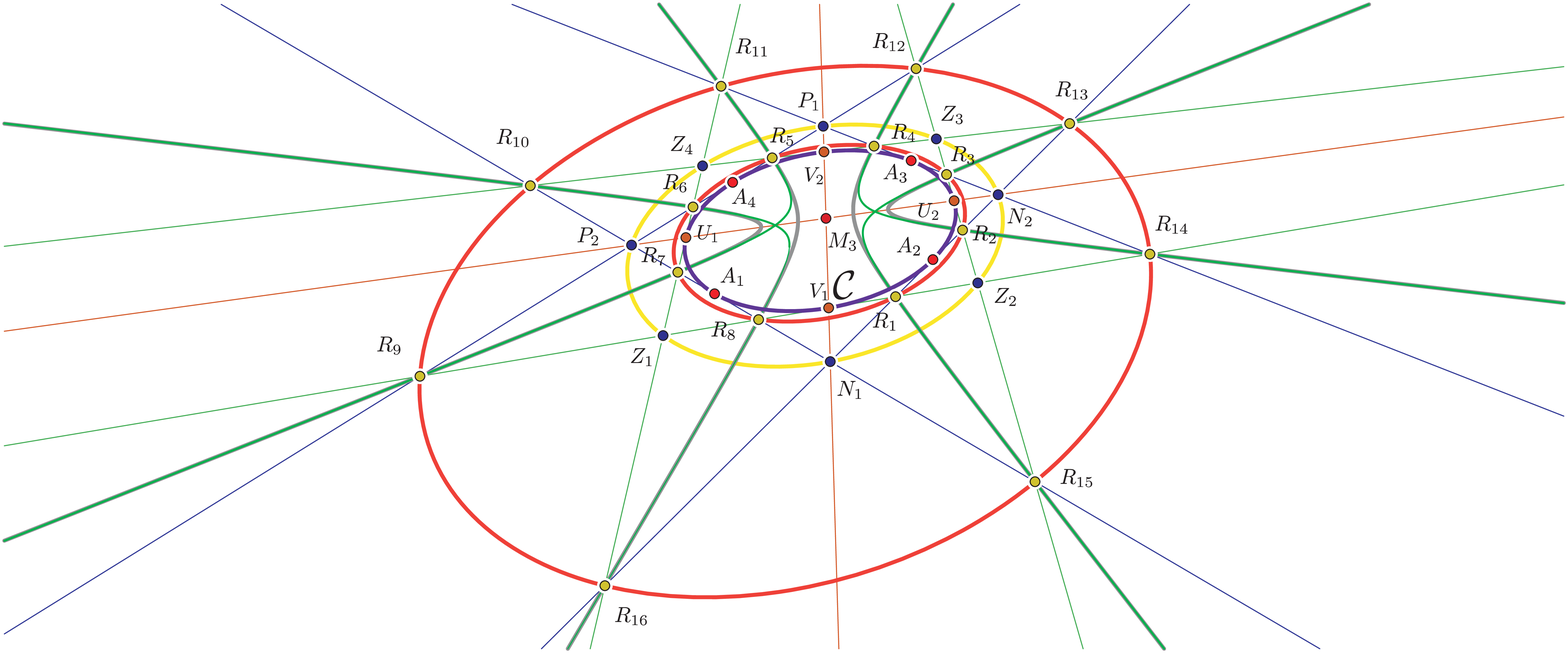, width=\textwidth}}
\caption{Theorem \ref{grid}} \label{next}
\end{figure}

Our theorems resemble Darboux's theorem, see \cite{Darbu}. They
could be seen as a very special case of Dragovi\'{c}-Radnovi\'{c}
theorem $8.38$, \cite{DraRad}. Such constructions are also studied
in the paper of Schwartz, see \cite{Schva1}. The following result
further explains their connection, but first we define $16$ points
of the intersections $R_1=l (Z_1 Z_2)\cap l(N_1 N_2)$, $R_2=l (Z_1
Z_2)\cap l(N_1 P_2)$, $R_3=l (Z_2 Z_3)\cap l(N_1 P_2)$, $R_4=l
(Z_2 Z_3)\cap l(P_1 P_2)$, $R_5=l (Z_3 Z_4)\cap l(P_1 P_2)$,
$R_6=l (Z_3 Z_4)\cap l(P_1 N_2)$, $R_7=l (Z_1 Z_4)\cap l(P_1
N_2)$, $R_8=l (Z_1 Z_4)\cap l(N_1 N_2)$, $R_9=l (Z_1 Z_2)\cap
l(P_1 P_2)$, $R_{10}=l (Z_3 Z_4)\cap l(N_1 P_2)$, $R_{11}=l (Z_2
Z_3)\cap l(P_1 N_2)$, $R_{12}=l (Z_1 Z_4)\cap l(P_1 P_2)$,
$R_{13}=l (Z_3 Z_4)\cap l(N_1 N_2)$, $R_{14}=l (Z_1 Z_2)\cap l(P_1
N_2)$, $R_{15}=l (Z_1 Z_4)\cap l(N_1 P_2)$ and $R_{16}=l (Z_2
Z_3)\cap l(N_1 N_2)$, see Figure \ref{next}.

\begin{theorem}\label{grid} The next groups of $8$ points lie on the same
conic:\\ $\{R_1, R_2, R_3, R_4, R_5, R_6, R_7, R_8\}$, $\{R_9,
R_{10}, R_{11}, R_{12}, R_{13}, R_{14}, R_{15}, R_{16}\}$,
\\ $\{R_1, R_2, R_5, R_6, R_{11}, R_{12}, R_{15}, R_{16}\}$, $\{R_3,
R_4, R_7, R_8, R_9, R_{10}, R_{13}, R_{14}\}$,\\ $\{R_1, R_5, R_7,
R_9, R_{11}, R_{13}, R_{15}\}$ and $\{R_2, R_3, R_4, R_6, R_8,
R_{10}, R_{12}, R_{14}, R_{16}\}$.
\end{theorem}

The proof of Theorem \ref{grid} uses the same arguments we used in
the previous proofs so we omit it.

If we look at the conic $\mathcal{C}$ and the conic $\mathcal{F}$
through the points $\{R_1, R_2, R_3, R_4, R_5$,\\ $R_6, R_7,
R_8\}$ we see they are Poncelet 8-connected and appropriate conics
from Theorem \ref{grid}, conic from Theorem \ref{tankon} with the
line $M_1 M_2$ form Poncelet-Darboux grid. Two conics $\{R_2, R_3,
R_4, R_6, R_8, R_{10}, R_{12}, R_{14}, R_{16}\}$ and $\{R_1, R_5,
R_7, R_9, R_{11}, R_{13}, R_{15}\}$ are not coming from
Poncelet-Darboux grid, but they could be directly obtained from
Dragovi\'{c}-Radnovi\'{c} theorem $8.38$, \cite{DraRad}. This
result improves the result of Schwartz \cite{Schva1} in a
particular case.

\section{Few words about Tabachnikov-Schwartz Theorem 4c}

In the end we think it is suitable to say something about already
mentioned Theorem 4c stated in \cite{Taba}. Tabachnikov and
Schwartz asked us for the proof. For this occasion we reformulate
it in the following manner:

\begin{theorem}[Tabachnikov-Schwartz Theorem 4c] Let $A_1 A_2
\dots A_{12}$ be a $12$-gon inscribed in a conic $\mathcal{C}$.
Let $\pi$ maps 12-gon $X_1 X_2 \dots X_{12}$ onto a new $12$-gon
according to the rule $\pi (X_i)=l (X_i X_{i+4})\cap l (X_{i+1}
X_{i+5})$. Then, 12-gon $A_1 A_2 \dots A_{12}$ is mapped with
$\pi^{(3)}=\pi\circ \pi\circ \pi$ onto a 12-gon inscribed in a
conic.
\end{theorem}

\begin{figure}[h!h!]
\centerline{\epsfig{figure=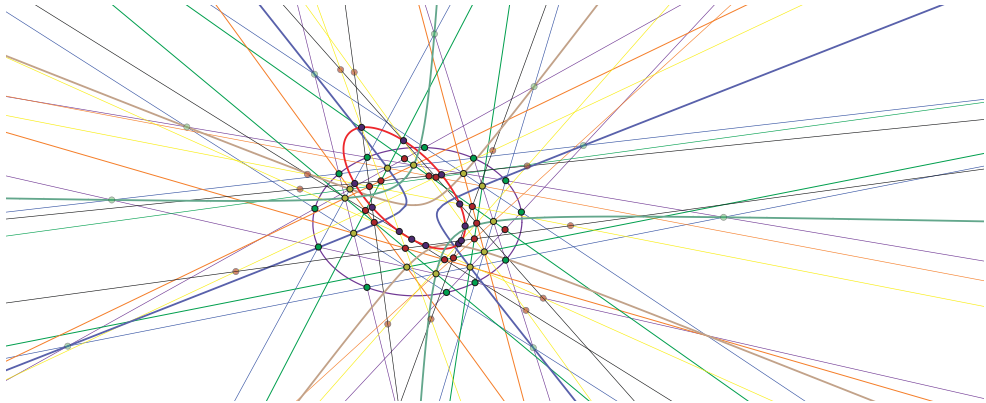,
width=\textwidth}} \caption{Tabachnikov-Schwartz Theorem 4c}
\label{tabaskica}
\end{figure}

This theorem was the starting point of our research. It seemed
that this theorem is a perfect candidate to use the technique
illustrated in \cite{BarIgi}, although in an unpublished paper of
Tabachnikov \cite{Taba1} one can find nice proofs for the theorems
from \cite{Taba}. Encouraged by our previous success, we tried to
prove Theorem 4c. We used \textit{Cinderella} again to test the
result and to obtain a nice picture. But at the beginning we
present the problem. We will explain Figure \ref{tabaskica}
carefully. We start with a 12-gon $A_1 A_2 \dots A_{12}$ (the
green points lying on the violet conic) inscribed in a conic and
define the (yellow) points obtained by $\pi$, (blue and violet
lines), $\pi^{(2)}$ the red points (green and orange lines) and
$\pi^{(3)}$ the violet points (black and yellow lines). It looks
like that at the every step we have a $6\times 6$ cage of curves,
see \cite{Katz}. But instead of dealing with 24 points at the
second step we take only $12$ of them. It is not possible to catch
the curves we want in the cage. By Mystic Octagon theorem we could
catch three interesting conics and one quartic in the blue-violet
cage. What to do with curves at other steps. Definitely we should
try to add some new points and then apply B\'{e}zout's theorem or
a similar statement. But what are that points and how to find
them? If we look more carefully, three quadrilaterals  can be
noticed ($A_1 A_4 A_7 A_{10}$, $A_2 A_5 A_8 A_{11}$ and $A_3 A_6
A_9 A_{12}$) inscribed in a conic  and usually the steps are
always defined as the certain intersection points of the side
lines of quadrilaterals. Thus, we thought if we want to overcome
the problems we faced, it is good to understand the quadrilaterals
in a conic better.

We have not succeeded in proving the Theorem 4c. But we conducted
some experiments in \textit{Cinderella} that we think are
important. Firstly, usually algebro-geometric facts give us some
freedom (for example, a product of $n$ lines could be often
generalized to a curve of degree $n$, see \cite{BarIgi}, etc.) but
here we have not found any such generalizations. Also, the
technique in \cite{BarIgi} usually does not differ order of
points, that means that certain permutations lead to new objects
of the same type (for example Pascal lines). Due to the difference
of three quadrilaterals we did not find new conic at the third
step. After all these experiments we believe Tabachnikov and
Schwartz Theorem 4c is more surprising and deeper fact then it
looks at the first glance!

\begin{center}\textmd{Acknowledgements }
\end{center}

The authors are grateful to Sergey Tabachnikov who brought Theorem
4c to our attention and sent us his unpublished reprints regarding
the subject. Special thanks to Vladimir Dragovi\'{c} and Milena
Radnovi\'{c} who helped us with discussions, comments and
encouragement.

\medskip

{\small \DJ{}OR\DJ{}E BARALI\'{C}, Mathematical Institute SASA,
Kneza Mihaila 36, p.p.\ 367, 11001 Belgrade, Serbia

E-mail address: djbaralic@mi.sanu.ac.rs

 BRANKO GRBI\'{C}, Mathematical Grammar School, Kraljice Natalije 37, 11000 Belgrade, Serbia

E-mail address: grbicbranko1@gmail.com

 \DJ{}OR\DJ{}E \v{Z}IKELI\'{C}, Mathematical Grammar School, Kraljice Natalije 37, 11000 Belgrade, Serbia

E-mail address: zikadj@hotmail.com}

\end{document}